\documentclass{mrlart3}
 \def\volno{xx}
  \def\yrno{2007}
  \def\issueno{x}
  \def\lpageno{100NN} 
\newtheorem{theorem}{Theorem}
\newtheorem*{thm}{Theorem}

\newtheorem{lem}[theorem]{Lemma}

\theoremstyle{remark}
\newtheorem{rem}[theorem]{Remark}

\newcommand{\wpi}{weighted Poincar\'e inequality}

\renewcommand{\bar}{\overline}

\newcommand{\eps}{\varepsilon}

\newcommand{\nonum}{\nonumber }
\newcommand{\setm}{\setminus }
\newcommand{\goto}{\rightarrow }

\newcommand{\re}{{\mathbb R}}

\newcommand{\bea}{\begin{eqnarray} }
\newcommand{\eea}{\end{eqnarray} }
\newcommand{\beay}{\begin{eqnarray*} }
\newcommand{\eeay}{\end{eqnarray*} }
\newcommand{\benu}{\begin{enumerate}}
\newcommand{\eenu}{\end{enumerate}}
\newcommand{\barr}{\begin{array}}
\newcommand{\earr}{\end{array}}
\newcommand{\grad }{\nabla}
\newcommand{\lap}{\triangle }
\newcommand{\pa }{\partial }
\newcommand{\sub }{\subseteq}

\newcommand{\p }{\mathcal{P}}
\newcommand{\lvl}{\mathcal {L}}

\newcommand{\la}{\langle}
\newcommand{\ra}{\rangle}

\newcommand{\lt}{\left}
\newcommand{\rt}{\right}
\newcommand{\bay}{\begin{array}}
\newcommand{\eay}{\end{array}}
\newcommand{\ric}{\mbox{Ric}}
\newcommand{\beg}{\begin}

\title{Finiteness and vanishing results on the weighted Poincar\'e inequality of complete manifolds}
\author{Kwan-hang Lam}

\begin{document}
\maketitle \author

\begin{abstract} We study manifolds satisfying a \wpi,
which was first introduced by Li-Wang \cite{liwangweight}. We generalized one of their
results by relaxing the Ricci curvature bound condition only being satisfied outside a
compact set and established a finitely many ends result. We also generalized a result of
\cite{liwangpos} to the \wpi\  case and established a vanishing result for $L^2$ harmonic
1-form provided that the weight function $\rho$ is of sub-quadratic growth of the
distance function.
\end{abstract}

\section*{Introduction}

In a work of Witten-Yau \cite{wy}, they proved that if $M^n$ is a conformally compact,
Einstein, $n\geq 3$ dimensional manifold whose boundary has positive Yamabe constant,
then $M$ must have only one end. Shortly after that, Cai-Galloway \cite{cg} relaxed the
assumption to allow the boundary of $M$ has non-negative Yamabe constant. In
\cite{xwang}, X. Wang proved the following:
\begin{thm}(X. Wang) Let $M^n$ be an n-dimensional ($n\geq 3$),
conformally compact manifold with Ricci curvature bounded from below by
$$Ric_M\geq -(n-1).$$ Let $\lambda_1(M)$ be the lower bound of the spectrum of
the Laplacian on $M$, If $$\lambda_1(M)\geq n-2,$$ then either\\
(a) $H^1(L^2(M))=0$; or\\
(b) $M=\mathbb R \times N$ with the warped product metric $ds^2=dt^2+\cosh ^2t\ ds_N^2,$
where $N$ is a compact manifold with $Ric_N\geq n-2.$ In particular, $M$ either has only
one end or it must be a warped product given as above.
\end{thm}
The above result generalized the results of Witten-Yau and Cai-Galloway since a theorem
of Mazzeo \cite{mazzeo}\ identifies the $L^2$ cohomology group $H^1(L^2(M))$ with the
relative cohomology group $H^1(M,\pa M)$ for conformally compact manifolds and a theorem
of Lee \cite{lee} asserts that $\lambda _1(M)=\frac{(n-1)^2}{4}$ for a $n$-dimensional
Einstein, conformally compact manifold with non-negative Yamabe constant for its
boundary. In \cite{liwangpos}, Li-Wang considered a complete, $n$-dimensional, Riemannian
manifold $M^n$ whose Ricci curvature is bounded from below by
$$Ric_M\geq -\frac{n-1}{n-2}\lambda_1(M),$$ where $\lambda _1(M)$, the greatest
lower bound of the spectrum of the Laplacian acting on $L^2$ functions, is assumed to be
positive. They proved the following
\begin{thm}(Li-Wang) Let $M^n$ be a complete Riemannian manifold of dimension $n\geq 3$. Suppose
$\lambda_1
(M)>0$ and $$Ric_M\geq -\frac{n-1}{n-2}\ \lambda_1 (M).$$ Then either \\
(1) $M$
has only one end with infinite volume; or\\
(2) $M=\re \times N$ with the warped product metric
$$ds_M^2=dt^2+\cosh^2\lt(\sqrt {\frac{\lambda _1(M)}{n-2}}\ t \rt)ds_N^2,$$ where
$N$ is a compact manifold with Ricci curvature bounded from below by
$$Ric_N\geq -\lambda_1 (M).$$
\end{thm}
Since all the ends of a conformally compact manifold must have infinite volume, the above
theorem thus generalized the work of \cite{cg}, \cite{xwang} and \cite{wy} to complete
manifolds with positive spectrum. Since $\lambda_1 (M)
>0,$ the variational principle for $\lambda_1 (M)$ implies the following
Poincar\'e inequality $$\lambda_1 (M)\int _M\phi ^2\leq \int_M|\nabla \phi |^2,$$ for any
compactly supported smooth function, $\phi \in C_c^\infty (M).$ In \cite{liwangweight},
the authors considered manifolds satisfying a weighted Poincar\'e inequality and
generalized most of their results in \cite{liwangpos} for manifolds with positive
spectrum to manifolds satisfying a weighted Poincar\'e inequality. A manifold $M^n$ is
said to be satisfied a weighted Poincar\'{e} inequality with a non-negative weight
function $\rho (x)$ if
$$\int _M\rho (x)\phi ^2(x)dV\leq \int _M|\nabla \phi
|^2dV, $$ for any compactly supported smooth function $\phi \in C_c^\infty (M).$ In
particular, when $\rho (x)=\lambda_1 (M)$ is a positive constant, $M$ is a manifold with
positive spectrum. We say that a  manifold $M$ has property $(\p_\rho)$ if a weighted
Poincar\'e inequality is valid on $M$ with some non-negative weight function $\rho $ and
the $\rho$-metric, defined by
$$ds_\rho ^2=\rho  \ ds_M^2$$ is complete. Define
$$S(R)=\sup_{B_\rho (R)}\sqrt \rho,$$ where $B_\rho(R)$ is the ball of radius
$R$ (with respect to some fixed point $p$) under the $\rho$-metric. In
\cite{liwangweight}, the authors proved the following:

\begin{theorem}\label{thm:liwang weighted_structural}(Li-Wang) Let $M^n$ be a complete manifold with dimension $n\geq 3$.
Assume that $M$ satisfies property $(\p_\rho )$ for some nonzero weight function $\rho
\geq 0$. Suppose $$Ric_M(x)\geq -\frac{n-1}{n-2}\ \rho (x)$$ for all $x\in M$. If $\rho $
satisfies the growth estimate
$$\liminf _{R\rightarrow \infty } \frac{S(R)}{F(R)}=0,
$$ where $$F(R)=\left\{\begin{array}{ll}\exp(\frac{n-3}{n-2}R)&\mbox{ for }n\geq 4\\R&\mbox{ for }n=3\end{array},\right.$$ then either
\begin{enumerate}\item $M$ has only one nonparabolic end; or
\item $M$ has two nonparabolic ends and is given by $M=\re \times
N$ with the warped product metric $$ds_M^2=dt^2+\eta ^2(t)ds_N^2,$$ for some positive
function $\eta (t)$, and some compact manifold $N$. Moreover, $\rho (t)$ is a function of
$t$ alone satisfying
$$\eta ''\eta ^{-1}=\rho $$ and $$\liminf _{x\goto \infty }\rho
(x)>0;\mbox{ or }$$
\item $M$ has one parabolic end and one nonparabolic end and is
given by $M=\re \times N$ with the warped product metric
$$ds_M^2=dt^2+\eta ^2(t)ds_N^2,$$ for some positive function $\eta (t)$, and some compact manifold
$N$. Moreover, $\rho (t)$ is a function of $t$ alone satisfying
$$\eta ''\eta ^{-1}=\rho $$ and $$\liminf _{x\goto \infty }\rho
(x)>0 \mbox{ on the nonparabolic end.}$$
\end{enumerate}
\end{theorem}
It is interesting to see if a similar theorem holds by relaxing the above assumptions to
be only satisfied outside a compact set of $M$. In this article, the following theorem
has been established:
\begin{theorem}\label{thm:lam weighted_finiteness}
Let $M^n$ be a complete manifold with dimension $n\geq 3$. Assume that $M$ has property
$(\p_\rho )$. Also assume that the Ricci curvature of $M$ satisfies the lower bound
$$Ric_{M\setminus K}(x)\geq -\frac{n-1}{n-2}\rho (x)+\tilde \varepsilon $$ for some $\tilde \varepsilon > 0$, compact set $K\subseteq M$. If $\rho $ satisfies the growth estimate
$$\liminf _{R\rightarrow \infty } \frac{S(R)}{F(R)}=0,
$$ where $$F(R)=\left\{\begin{array}{ll}\exp(\frac{n-3}{n-2}R)&\mbox{ for }n\geq 4\\R&\mbox{ for }n=3\end{array},\right.$$
then $M$ has only finitely many non-parabolic ends. If in addition $\displaystyle \liminf
_{x\goto \infty}\rho (x)>0$, then $M$ has only finitely many infinite volume ends.
\end{theorem}

On the other hand, Li-Wang also proved a vanishing theorem for $L^2$ integrable harmonic
1-forms on $M$ in \cite{liwangpos}:

\beg{theorem}\label{thm:liwang pos_vanishing}(Li-Wang) Let $M$ be a $n$-dimensional
complete Riemannian manifold with $\lambda _1(M)>0$ and $$Ric_M\geq
-\frac{n}{n-1}\lambda_1 (M)+\eps ,$$ for some $\eps
>0.$ Then $H^1(L^2(M))=0.$
\end{theorem}

In this article, we generalized above theorem to manifolds satisfying a weighted
Poincar\'{e} inequality:

\begin{theorem}\label{thm:weighted1form} Let $M$ be a $n$-dimensional complete Riemannian manifold
satisfying a weighted Poincar\'e inequality with a non-negative weight function $\rho
(x).$ Assume the Ricci curvature satisfies $$Ric_M(x)\geq -\frac{n}{n-1}\ \rho(x) +\eps
,$$ for some $\eps >0.$ Let $r(x,p)$ be the distance function from $x$ to some fixed
point $p$, If $\rho (x)=O(r^{2-\alpha }(x,p)),$ for some $0<\alpha <2.$ Then
$H^1(L^2(M))=0.$
\end{theorem}

\section*{ Some lemmas}

The following lemma is modified from \cite{liwangweight} to suit our situation.
\begin{lem} \label{lem:grad}Assume that $$Ric _M(x)\geq -\frac{n-1}{n-2}\ l(x)$$
for some function $l(x)$. If $f$ is a positive harmonic function on $M$, then $$|\nabla
f|(x)\leq \left( (n-1)\sup _{B_\rho (x, 1)}\sqrt {l(y)}+C_1\sup_{B_\rho (x,1)}\sqrt {\rho
(y)}\right)f(x),$$ where $C_1$ is a constant only depending on $n$. In particular, if the
lower bound of the Ricci curvature of $M$ satisfies
$$Ric_{M\setminus K}(x)\geq -\frac{n-1}{n-2}\ \rho (x)+\tilde \varepsilon $$
where $K$ is a compact sub-domain of $M$ and for some $\tilde \varepsilon > 0$. Then\bea
\label{e:grad}|\grad f|(x)&\leq &C (\sup_{B_\rho (x,1)}\sqrt \rho )f(x), \eea where
$C=C(n)$, provided that $B_\rho (x,1)\cap K=\phi $ .\end{lem}
\begin{proof} Cheng-Yau's \cite{CY75} (see also \cite{liwangpos}) local gradient estimate for positive
harmonic functions implies that for any $R>0$, \bea \label{e:1a}|\grad f|(x)\leq \left(
(n-1)\sup _{B(x,R)}\sqrt {l }+CR^{-1}\right)f(x),\eea where $C=C(n)$. Consider the
function $g(r)=r-{(\sup _{B(x,r)}\sqrt \rho)}^{-1}$. Since $g$ is negative when $r\goto
0$ and $g\goto +\infty $ as $r\goto +\infty $. Hence we can choose $R_0>0$ such that
$g(R_0)=0$, that is, $R_0={(\sup _{B(x, R_0)}\sqrt \rho )}^{-1}$. For any point $y\in
B(x, R_0)$, let $\gamma $ be a minimizing geodesic (with respect to $ds_M^2$)
joining $x,y$, then \beay r_\rho (x,y)&\leq &\int_\gamma \sqrt {\rho (\gamma (t))}dt\\
&\leq &\left(\sup _{B(x,R_0)}\sqrt {(\rho (y))}\right)R_0\\
&\leq &1,\eeay hence $B(x,R_0)\subseteq B_\rho (x,1)$. Combining with (\ref{e:1a}) by
choosing $R=R_0$, the result follows.
\end{proof}
The following version of Bochner formula is well-known and was first used by Yau
\cite{Yau75}.
\begin{lem}\label{bochner}\cite{Yau75} (see also \cite{liwangweight}) Let
$M^n$ be a complete Riemannian manifold of dimension $n\geq 2$. Assume that the Ricci
curvature of $M$ satisfies the lower bound
$$Ric_M(x)\geq -(n-1)\tau (x).$$ Assume $f$ is a nonconstant
harmonic function on $M$. Then the function $h=|\grad f|$ satisfies
$$\lap h\geq -(n-1)\tau h+\frac{|\grad h|^2}{(n-1)h}.$$ In addition,
if we write $g=h^\frac{n-2}{n-1}$, the above inequality becomes
$$\lap g\geq -(n-2)\tau g.$$
\end{lem}
\begin{proof}For the sake of completeness, we outline the proof here. We choose a local orthonormal frame $\{e_1,e_2,\cdots , e_n\}$ such
that $e_1f=|\grad f|$ and $e_\alpha f =0$, for $\alpha =2,\cdots , n$ at
a point $x$. \beay h^2&=&|\grad f|^2\\
hh_j&=&\sum _{i=1}^nf_{ij}f_i. \eeay Hence \bea h_j&=&f_{1j},\nonum
\\\label{e:b1}|\grad h|^2&=&\sum_{j=1}^nh_j^2=\sum_{j=1}^nf_{1j}^2 , \eea when evaluate at
$x$. Using the fact that $f$ is harmonic and the Ricci formula, we compute \bea |\grad
h|^2+h\lap h&=&\sum_{j=1}^n(h_j^2+hh_{jj})\label{e:b2}
\\&=&\sum_{i,j=1}^n(f_{ij}^2+f_{ijj}f_i)\nonum \\
&=&\sum_{i,j=1}^n(f_{ij}^2+f_{jji}f_i+R_{ij}f_if_j)\nonum \\
&=&\sum_{i,j=1}^n(f_{ij}^2+R_{ij}f_if_j)\nonum \\
&\geq &-(n-1)\tau h^2+\sum_{i,j=1}^nf_{ij}^2, \nonum \eea and \beay \sum
_{i,j=1}^nf_{ij}^2&=&2\sum _{\alpha =2}^nf_{1\alpha
}^2+f_{11}^2+\sum_{\beta =2}^nf_{\beta \beta }^2\\
&\geq &2\sum _{\alpha =2}^nf_{1\alpha }^2+f_{11}^2+\frac{1}{n-1}\left(\sum_{\beta
=2}^nf_{\beta \beta
}\right)^2\\
&=&2\sum _{\alpha =2}^nf_{1\alpha }^2+\frac{n}{n-1}f_{11}^2\\
&\geq &\frac{n}{n-1}\sum _{j =1}^nf_{1j }^2.\eeay Combining the above inequality with
(\ref{e:b1}), (\ref{e:b2}) and evaluate at $x$, the result follows.
\end{proof}
A theory of Li-Tam \cite{litamstructure} allows us to count the number of non-parabolic
ends of $M$ by counting the dimension of $\mathcal K^0(M)$, a subspace of the space of
all harmonic functions on $M$. We outline the construction of Li-Tam here. Assume that
$M$ has at least two non-parabolic ends, $E_1, E_2$, for $R>0$, we solve the following
equation
$$\left\{\begin{array}{rlc}\lap
f_R=0&\mbox{ in }&B(R)\\
f_R=1&\mbox{ on }&\partial B(R)\cap E_1\\
f_R=0&\mbox{ on }&\partial B(R)\setm E_1.\end{array}\right.$$ By passing to a convergent
subsequence, the sequence $\{f_R\}$ converges to a nonconstant harmonic function $f_1$ with
finite Dirichlet integral, satisfying $0\leq f_1\leq 1$. Clearly for each non-parabolic end
$E_i$, we can construct a corresponding $f_i$ by the above process. Let $\mathcal K^0(M)$ be
the linear space containing all the $f_i$'s constructed as above. By the construction, the
number of nonparabolic ends of $M$ is given by the dimension of $\mathcal K^0(M)$. The
following lemma of Li is useful in proving finiteness type theorems.

\begin{lem}\cite{lidimensionlemma}\label{lem:dimension} Let $\mathcal{H}$ be a finite dimensional
subspace of $L^2$ $p$-forms defined over a set $D\sub M^n$. If $V(D)$ denotes the volume
of the set $D$, then there exists $\omega_0\in \mathcal{H}$ such that \beay \dim
\mathcal{H}\ \int_D|\omega_0|^2\leq V(D)\cdot \sup _D|\omega_0|^2\cdot
\min\left\{{n\choose p},\dim \mathcal{H}\right\}.\eeay
\end{lem}


\section*{Result of finitely many ends}


\begin{theorem}\label{thm:main1}
Let $M^n$ be a complete manifold with dimension $n\geq 3$. Assume that $M$ has property
$(\p_\rho )$. Also assume that the Ricci curvature of $M$ satisfies the lower bound
$$Ric_{M\setminus K}(x)\geq -\frac{n-1}{n-2}\ \rho (x)+\tilde \varepsilon $$ for some $\tilde \varepsilon > 0$,
compact set $K\subseteq M$. If $\rho $ satisfies the growth estimate
$$\liminf _{R\rightarrow \infty } \frac{S(R)}{F(R)}=0,
$$ where $$F(R)=\left\{\begin{array}{ll}\exp(\frac{n-3}{n-2}R)&\mbox{ for }n\geq 4\\R&\mbox{ for }n=3\end{array},
\right.$$
then $M$ has only finitely many nonparabolic ends. If in addition $\displaystyle \liminf
_{x\goto \infty}\rho (x)>0$, then $M$ has only finitely many infinite volume ends.
\end{theorem}

\begin{proof} By the discussion above lemma \ref{lem:dimension}, it is sufficient to estimate
$\dim \mathcal K^0(M)$. We may assume that $M$ has at least two nonparabolic ends and
hence there exists a nonconstant bounded harmonic function $f\in \mathcal K^0(M)$ with
finite Dirichlet integral. By maximal principle, we may assume that $\inf f=0$ and $\sup
f=1$ (see \cite{litamstructure}). Bochner formula and the assumption on the Ricci
curvature give us
$$\lap |\grad f|\geq -\frac{n-1}{n-2}|\grad f|(\rho -\eps) +\frac{|\grad |\grad f||^2}{(n-1)|\grad
f|},$$ where $\eps =\frac{n-2}{n-1}\tilde \eps$. Applying lemma \ref{bochner}, the above
equation becomes \bea \label{e:1}\triangle g+\rho g\geq \varepsilon g \eea in $M\setminus
K$, where $g=|\grad f|^\frac{n-2}{n-1}$. Let $\phi \in C_c^\infty (M\setm K)$ be a
non-negative smooth function with compact support in $M\setm K$. Using the property $(\p
_\rho )$ of $M$ and integration by parts, we have \bea \label{e:2}\int_M\phi ^2\rho
g^2&\leq &\int _M|\grad (\phi g)|^2\\&=&\int_M|\grad
\phi|^2g^2+\int_M\phi ^2|\grad g|^2+2\int_M\phi g\langle\grad \phi ,\grad g\rangle \nonum\\
&=&\int_M|\grad \phi|^2g^2+\int_M\phi ^2|\grad g|^2+\frac{1}{2}\int_M\langle \grad (\phi
)^2,\grad(g)^2\rangle \nonumber
\\
&=&\int_M|\grad \phi|^2g^2-\int_M\phi ^2g\lap g.\nonumber\eea Combining (\ref{e:1}) and
(\ref{e:2}), we have \bea \nonum \eps \int_M\phi ^2g^2\leq \int_M|\grad \phi |^2g^2,
\forall \phi \in C_c^\infty (M\setm K).\eea Since $K$ is compact, we may choose $R_0>0$
such that $$K\sub\bigcup _{x\in K}B_\rho (x,1)\subseteq B(R_0-1).$$ Let $R>0$ be such
that $B(R_0)\subseteq B_{\rho }(R-1)$. The above inequality implies \bea \label{e:3}\eps
\int_{B_{\rho (R)}\setm B(R_0-1)}\phi ^2g^2\leq \int_{B_\rho (R)\setm B(R_0-1)}|\grad
\phi |^2g^2, \eea for any $\phi \in C_c^\infty (B_\rho (R)\setm B(R_0-1))$. If we write
$\phi =\psi \cdot \chi $, the right hand side of (\ref{e:3}) becomes \bea \label{e:4}
\int_{B_\rho (R)\setm B(R_0-1)}|\grad \phi |^2g^2&\leq &2\int _{B_\rho (R)\setm
B(R_0-1)}|\grad \psi |^2\chi ^2g^2\\&+&2\int_{B_\rho (R)\setm B(R_0-1)} |\grad \chi
|^2\psi ^2g^2. \nonum \eea Now we choose $\chi ,\psi $ as in \cite{liwangweight}:
$$\chi (x)=\left\{\begin{array}{cl}0& \mbox{ on }\lvl(0,\delta \bar \eps )\cup \lvl(1-\delta \bar \eps ,1) \\
(-\log \delta )^{-1}(\log f-\log (\delta \bar \eps ))&\mbox{ on  } \lvl(\delta \bar \eps
,\bar \eps ) \cap (M\setm
E_1)\\
(-\log \delta )^{-1}(\log (1-f)-\log (\delta \bar \eps ))&\mbox{ on
}\lvl(1-\bar \eps ,1-\delta \bar \eps )\cap E_1\\
1&\mbox{ otherwise}\end{array}\right. ,$$ for some $0<\delta <1$ and $0<\bar \eps
<\frac{1}{2}$ to be determined later, where \beay \lvl(a,b)=\{x\in M:a<f(x)<b\}.\eeay
$$\psi
(x)=\left\{\begin{array}{cl}0&\mbox{ on }B(R_0-1)\\1&\mbox{ on
}B_\rho (R-1)\setm B(R_0)\\
R-r_\rho &\mbox{ on } B_\rho (R)\setm B_\rho(R-1)\\
0&\mbox{ on }M\setm B_\rho (R)
\end{array}.\right.$$ The first term on the right hand side of
(\ref{e:4}) can be estimated by \bea \label{e:5}2\int _{B_\rho (R)\setm B(R_0-1)}|\grad
\psi |^2\chi ^2g^2&\leq & 2\int _{B_\rho (R)\setm B_\rho (R-1)}|\grad \psi |^2\chi
^2g^2\\&+&2\int _{B(R_0)\setm B(R_0-1)}|\grad \psi |^2\chi ^2g^2\nonum .\eea Now consider
\bea &&\ \ \ \ \ \ \int_{E_1}|\grad \chi
|^2\psi ^2g^2\\
&\leq &(\log \delta )^{-2}\int _{\lvl(1-\bar \eps, 1-\delta \bar \eps )\cap E_1\cap
(B_\rho (R)\setm B(R_0-1))}|\grad
f|^{2+\frac{2(n-1)}{n-2}}(1-f)^{-2}\nonum \\
&\leq &C\ S^{\frac{2(n-2)}{n-1}}(R+1)(\log \delta )^{-2}\int _{\lvl(1-\bar \eps, 1-\delta
\bar \eps )\cap E_1\cap (B_\rho (R)\setm
B(R_0-1))}|\grad f|^2(1-f)^{\frac{2(n-2)}{n-1}-2}\nonum \\
&\leq &C\ S^{\frac{2(n-2)}{n-1}}(R+1)(\log \delta )^{-2}\int _{\lvl(1-\bar \eps, 1-\delta
\bar \eps )\cap E_1\cap B_\rho (R)}|\grad f|^2(1-f)^{\frac{2(n-2)}{n-1}-2}\nonum ,\eea
where the second inequality follows from (\ref{e:grad}) by replacing $f$ by $1-f$ and the
choice of $R_0$ so that $B_\rho (x,1)$ does not intersect $K$ for any $x\in B_\rho
(R)\setminus B(R_0-1).$ Now, we are exactly the same situation as of \cite{liwangweight}.
We follow the choices of $\delta , \bar \eps $ as in \cite{liwangweight}, $\bar \eps
\goto 0$ as $R\goto +\infty$ and we have the followings: \beay \int_M|\grad \chi |^2\psi
^2g^2 \goto 0 &\mbox{ as
}& R\goto +\infty , \mbox{ and }\\
\int _{B_\rho (R)\setm B_\rho (R-1)}|\grad \psi |^2\chi ^2g^2\goto 0&\mbox{ as }&R\goto
+\infty .\eeay Now if we let $R\goto +\infty $ and combining the above results with
(\ref{e:4}) and (\ref{e:5}), (\ref{e:3}) becomes \bea \nonum \eps \int _{M\setm
B(R_0)}g^2\leq \tilde C \int_{B(R_0)\setm B(R_0-1)}g^2, \eea where $\tilde C$ is a
constant depending only on $n$, which implies \bea \label{e:7}\int _{B(2R_0)}g^2\leq C
\int_{B(R_0)}g^2,\eea where $C=C(\eps , n)$. Since the function $g$ satisfies the
differential inequality
$$\lap g\geq -\alpha g$$ on $B(p, 2R_0)$, where $\alpha =\inf
_{B(p,2R_0)}Ric_M$, the mean value inequality of Li-Tam \cite{LT91} gives
\beay g^2(x)&\leq &C_1\int_{B(x, R_0)}g^2\\
&\leq &C_1 \int _{B(p,2R_0)}g^2\eeay for any $x\in B(p,R_0)$, where $\nu =\inf _{x\in
B(p,R_0)}V_x(R_0)$ and $C_1=C_1(n, \alpha , \nu )$. Combining with (\ref{e:7}), we have
$$\sup _{B(p, R_0)}g^2\leq C_2\int _{B(p, R_0)}g^2,$$ where $C_2=C_2(\eps, n,\alpha, \nu).$ On the other hand, the
Schwarz's inequality implies that $$\int _{B(R_0)}g^2\leq \left(\int_{B(R_0)}|\grad
f|^2\right)^{\frac{n-2}{n-1}}V_p(R_0)^{\frac{1}{n-1}}.$$ Therefore we have \bea
\label{e:8}\sup _{B(R_0)}|\grad f|^2\leq C_3\int _{B(R_0)}|\grad f|^2,\eea where
$C_3=C_3(\eps, n,\alpha, \nu, R_0)$ is a constant independent of $f\in \mathcal K^0(M)$.
By unique continuation, $$\int _{B(R_0)}|\grad f|^2\not=0,$$ provided that $f$ is not a
constant function. Therefore,
$$\int_{B(R_0)}\la \grad f, \grad g\ra$$
defines a non-degenerate bilinear form on the space of $1$-forms $$\mathcal K=\{df:f\in
\mathcal K^0(M)\}.$$ Lemma \ref{lem:dimension} asserts that there exists $f_0\in \mathcal
K^0(M)$ such that
$$\dim \mathcal K\int _{B(R_0)}|df_0|^2\leq n V_p(R_0)\sup _{B(R_0)}|df_0|^2.$$
Combining the above with (\ref{e:8}) implies
$$\dim \mathcal K^0(M)=\dim \mathcal K+1\leq C_4$$
for some fixed constant $C_4=C_4(C_3,V_p(R_0))$, which completes the proof. The second
part of the theorem follows from \cite{liwangweight}, an end is nonparabolic if and only
if it has infinite volume, provided that $\displaystyle \liminf _{x\goto \infty}\rho
>0$.
\end{proof}
\begin{rem} We would like to point out that the growth condition is not too
restrictive. It is satisfied if the weight function $\rho $ does not growth too fast, for
instance, if $\rho $ has only polynomial growth ($n\geq 4$).
\end{rem}


\section*{Vanishing theorems on $L^2$ harmonic forms} In this section we study $H^1(L^2(M)),$
the space of $L^2$ integrable harmonic 1-forms. If $f$ is a harmonic function with finite
Dirichlet integral, then the exterior derivative of $f$, $df$ is a $L^2$ integrable harmonic
1-form. By the theory of Li-Tam \cite{litamstructure} (see also \cite{liwangpos}), we
have \beay \dim H^1(L^2(M))+1&\geq &\dim \mathcal {K}^0(M)\\
&\geq &\mbox{number of non-parabolic ends of $M$}.\eeay If we further assume that
$\lambda_1(M)>0,$ then
$$\dim H^1(L^2(M))+1\geq \mbox{number of infinite volume ends of $M$}.$$
Therefore, an estimate on $\dim H^1(L^2(M))$ is, in general, a stronger estimate than an
estimate on the number of nonparabolic ends (infinite volume ends if $\lambda_1(M)>0$). It is
known that if $\omega $ is a $L^2$ harmonic 1-form, then it is both closed and co-closed. In
particular, $h=|\omega |$ satisfies a Bochner type formula
$$\lap h\geq \frac{\mbox{Ric}_M(\omega, \omega )}{h}+\frac{|\grad
h|^2}{(n-1)h}.$$ We start by proving an estimate for functions that satisfy the above Bochner
type formula. We believe the estimate is of independent interest and will be useful in many
other situations.
\begin{lem}\label{lem:est} Let $b >-1$. Assume $h$ satisfies differential inequality $$\lap h\geq -ah+b\frac{|\grad
h|^2}{h},$$ for some constant $a$. For any $\eps
>0$, we have the following
estimate $$\left(b(1-\eps )+1\right)\int _M|\grad (\phi h)|^2\leq \left(b(\frac{1}{\eps
}-1)+1\right)\int _Mh^2|\grad \phi |^2+a\int_M\phi ^2h^2,
$$ for any compactly supported smooth function $\phi \in C_c^\infty (M)$. In addition,
if $$\int _{B_p(R)}h^2=o(R^2),$$ then
$$\int_M|\grad h|^2\leq \frac{a}{b+1}\int_Mh^2 .$$ In particular, $h$ has finite Dirichlet
integral if $h\in L^2(M).$
\end{lem}

\begin{proof} Let $\phi \in C_c^\infty (M)$ be a smooth function
with compact support. Integration by parts implies \bea \label{e:est1}\int_M\phi ^2h\lap
h&=&-\int_M\langle \grad (\phi
^2h),\grad h\rangle\\
\nonum&=&\int_M\langle h\grad \phi +\grad (\phi h),h\grad \phi
-\grad (\phi h)\rangle\\
\nonum &=&\int_Mh^2|\grad \phi |^2-\int_M|\grad (\phi h)|^2.\eea By the differential
inequality satisfied by $h$, we have
$$\int_M\phi ^2h\lap h\geq -a\int_M\phi ^2 h^2+b\int _M\phi ^2|\grad h|^2.$$ Combining
the above inequality with (\ref{e:est1}) gives \bea \label{e:est2}\int _Mh^2|\grad \phi
|^2+a\int_M\phi ^2h^2\geq b\int _M\phi ^2|\grad h|^2+\int_M|\grad (\phi h)|^2.\eea On the
other hand, Schwarz inequality implies \beay \int_M\phi ^2|\grad h|^2&=&\int_M\langle \grad
(\phi h)-h\grad \phi ,
\grad (\phi h)-h\grad \phi \rangle \\
&=& \int_M(|\grad (\phi h)|^2-2\langle \grad (\phi h),h\grad \phi\rangle+h^2|\grad \phi |^2)\\
&\geq &(1-\eps )\int_M|\grad (\phi h)|^2+(1-\frac{1}{\eps })\int_Mh^2|\grad \phi |^2, \eeay
for any $\eps >0.$ Combining the above with (\ref{e:est2}) gives the first result of the
lemma. For the second part, we choose $$\phi =\left\{\begin{array}{ll}1&\mbox{ on
}B(R)\\
0 &\mbox{ on }M\setm B(2R)\end{array}\right.$$ such that $|\grad \phi |\leq C/R$ on
$B(2R)\setm B(R)$. The estimate of the lemma implies $$\left(b(1-\eps )+1\right)\int
_{B(R)}|\grad h|^2\leq \left(b(\frac{1}{\eps }-1)+1\right)R^{-2}\int _{B(2R)\setminus
B(R)}h^2+a\int_Mh^2.$$ Using $\int_{B_p(R)}h^2=o(R^2),$ and let $R\goto +\infty, \eps \goto
0$, the second part of the lemma is achieved.

\end{proof}

\begin{rem} By the proof of the above lemma, the above conclusions are still valid if we
assume $a=a(x)$ is a function of $x$, in the following forms
$$\left(b(1-\eps )+1\right)\int _M|\grad (\phi h)|^2\leq \left(b(\frac{1}{\eps }-1)+1\right)\int _Mh^2|\grad \phi |^2+\int_Ma(x)\phi
^2h^2 ;$$ and $$\int_M|\grad h|^2\leq \frac{1}{b+1}\int_Ma(x)h^2 .$$
\end{rem} The following theorem is an immediate application of the above lemma.

\beg{cor}\label{cor:vanish} With all the assumptions in lemma \ref{lem:est}, if
$\lambda_1(M)>0$ and the Ricci curvature satisfies
$$\ric_M\geq -(b+1)\lambda_1(M)+\delta,$$ for some $\delta
>0$. Then $H^1(L^2(M))=0.$ In particular, we recover a theorem of Li-Wang in \cite{liwangpos}: Let $M^n$ be a complete noncompact $n$-dimensional Riemannian manifold with
$\lambda_1(M)>0,$ and the Ricci curvature satisfies
$$\ric_M\geq -\frac{n}{n-1}\lambda_1(M)+\delta,$$
for some $\delta >0.$ Then $H^1(L^2(M))=0.$\end{cor} \beg{proof} Let $\omega \in H^1(L^2(M))$
and $h=|\omega |\in L^2(M)$. Bochner formula implies $h$ satisfies the differential inequality
$$\lap h\geq -ah+b\frac{|\grad h|^2}{h},$$ with $a=(b+1)\lambda_1(M)-\delta
$. Combining lemma \ref{lem:est} with the variational principle of $\lambda_1(M)$, we have
\beay (b(1-\eps) +1))\lambda_1(M)\int_M \phi
^2h^2 &\leq &(b(1-\eps) +1))\int_M|(\grad \phi h)|^2\\
&\leq &a\int_M\phi ^2h^2+\left(b(\frac{1}{\eps }-1)+1\right)\int _Mh^2|\grad \phi |^2,\eeay
for any $\eps
>0$ and any compactly supported smooth function $\phi \in C_c^\infty (M).$ The above inequality implies \beay \delta \int_M\phi ^2h^2\leq b\eps
\lambda_1(M)\int_M\phi^2h^2+\left(b(\frac{1}{\eps }-1)+1\right)\int _Mh^2|\grad \phi |^2.\eeay
Let
$$\phi =\left\{\begin{array}{ll}1&\mbox{ on
}B(R)\\
0 &\mbox{ on }M\setm B(2R)\end{array}\right.$$ such that $|\grad \phi |^2\leq C/R^2$ on
$B(2R)\setm B(R)$. The above inequality becomes $$\delta \int_{B(R)}h^2\leq b\eps
\lambda_1(M)\int_{B(2R)}h^2+R^{-2}\left(b(\frac{1}{\eps }-1)+1\right)\int_{B(2R)\setminus
B(R)}h^2.$$ Combining the above inequality with the assumption $h\in L^2(M)$ and let $R\goto
+\infty $ and then take $\eps \goto 0$, we conclude that $\int_Mh^2\leq 0$ and thus $h\equiv
0.$ For the second part, we just need to notice that in general, for a Riemannian manifold
$M^n$ of dimension $n$, Bochner formula for $h=|\omega |$ is valid with $b=\frac{1}{n-1}.$
\end{proof}

\begin{theorem}\label{thm:weighted1form} Let $M^n$ be a complete noncompact
manifold of dimension $n$ satisfying the weighted Poincar\'e inequality with a non-negative
weight function $\rho (x).$ Assume the Ricci curvature satisfies
$$Ric_M(x)\geq -\frac{n}{n-1}\rho(x) +\delta ,$$ for some $\delta >0.$ If
$\rho (x)=O(r_p^{2-\alpha }(x)),$ where $r_p(x)$ is the distance function from $x$ to some
fixed point $p$, for some $0<\alpha <2$. Then $H^1(L^2(M))=0.$
\end{theorem}
\begin{proof} Let $\omega \in H^1(L^2(M))$
and $h=|\omega |\in L^2(M)$. Applying Lemma \ref{lem:est} with $b=\frac{1}{n-1}, \ a=(b+1)\rho
-\delta$ and using the \wpi\  we have \beay (b(1-\eps) +1))\int_M\rho \phi
^2h^2 &\leq &(b(1-\eps) +1))\int_M|(\grad \phi h)|^2\\
&\leq &\int_Ma(x)\phi ^2h^2+\left(b(\frac{1}{\eps }-1)+1\right)\int _Mh^2|\grad \phi |^2.\eeay
Thus
$$\delta \int_M\phi^2h^2\leq \left((1-b)+b\eps ^{-1}\right)\int_Mh^2|\grad \phi |^2+b\eps \int_M\rho \phi^2h^2,$$ for any $\eps
>0$ and any compactly supported smooth function $\phi \in C_c^\infty (M).$ Let $\eps =R^{\alpha /2-2}$ and $$\phi
=\left\{\begin{array}{ll}1&\mbox{ on
}B(R)\\
0 &\mbox{ on }M\setm B(2R)\end{array}\right.$$ such that $|\grad \phi |^2\leq C/R^2$ on
$B(2R)\setm B(R)$. The above inequality becomes $$\delta \int_{B(R)}h^2\leq bR^{-\alpha
/2}\int_{B(2R)}h^2+\left((1-b)R^{-2}+bR^{-\alpha /2}\right)\int_{B(2R)}h^2.$$ Combining the
above inequality with the assumption $h\in L^2(M)$ and let $R\goto +\infty $, we conclude that
$\int_Mh^2\leq 0$ and thus $h\equiv 0.$
\end{proof}
\begin{theorem}\label{thm:weighted1form2} Let $M^n$ be a complete noncompact
manifold of dimension $n$ satisfying the weighted Poincar\'e inequality with weight function
$\rho >0.$ Assume the Ricci curvature satisfies $Ric_M\geq -\left(\frac{n}{n-1}-\delta
\right)\rho ,$ for some $\delta >0,$ and $\rho =O(r_p^{2-\alpha }(x))$ for some $0<\alpha <2$.
Then $H^1(L^2(M))=0.$
\end{theorem}
\begin{proof} Let $\omega \in H^1(L^2(M))$ and $h=|\omega |\in L^2(M)$.
Applying Lemma \ref{lem:est} with $b=\frac{1}{n-1},\ a=(b+1-\delta )\rho,$ and \wpi we have
\beay (b(1-\eps) +1))\int_M\rho \phi ^2h^2 &\leq &(b(1-\eps) +1))\int_M|(\grad \phi
h)|^2\\
&\leq &a\int_M\phi ^2h^2+\left(b(\frac{1}{\eps }-1)+1\right)\int _Mh^2|\grad \phi |^2.\eeay
Thus $$\delta \int_M\rho \phi^2h^2\leq b\eps \int_M\rho
\phi^2h^2+(b\eps^{-1}-b+1)\int_Mh^2|\grad \phi|^2.$$ If we choose
$$\phi =\left\{\begin{array}{ll}1&\mbox{ on
}B(R)\\
0 &\mbox{ on }M\setm B(2R)\end{array}\right.$$ Using $\rho(x)=O(r^{2-\alpha }(x,p))$ and let
$\eps =R^{\alpha /2-2}$, the above inequality implies
$$\delta \int_{B(R)}\rho h^2\leq CR^{-\alpha /2}\int_{B(2R)}h^2+(bR ^{-\alpha
/2}+(1-b)R^{-2})\int_{B(2R)}h^2,$$for some constant $C$. Using $h\in L^2(M)$ and letting
$R\goto +\infty$, we conclude that $\int_M\rho h^2=0.$ By Lemma \ref{lem:est}, we have \beay
\int _M|\grad h|^2&\leq &\left(1-\frac{\delta }{b+1}\right)\int_M\rho
h^2\\
&= &0.\eeay Therefore $|\grad h|\equiv 0$. Hence $h\equiv c\in L^2(M)$ for some constant $c$.
Since $M$ is non-parabolic, it must have infinite volume and thus $h\equiv 0.$
\end{proof}


  \bibliographystyle{mrl}
  \bibliography{bib}

\noindent
\address{Department of Mathematics \\ National Cheng Kung University, Taiwan \\ 1 university
road, Tainan 701, Taiwan}

\noindent \email{khlam@alumni.uci.edu}

  \end{document}